\documentclass[11 pt, psamsfonts]{amsart}
\usepackage{amssymb}

\usepackage{myart}

\def\um{U^-}
\def\umq{U^-_\Q}
\def\usl{U_v(\widehat{\mathfrak {sl}}_l)}
\def\g{\mathfrak g}
\def\wt{{\rm{wt}}}

\def\R{{\mathbb R}}
\def\Z{{\mathbb Z}}
\def\nat{{\mathbb N}}
\def\Q{{\mathbb Q}}
\def\C{{\mathbb C}}

\def\dlm{d_{\la, \mu }}
\def\alm{\alpha_{\la, \mu }}
\def\W{\mathcal W}
\def\S{\mathcal S}
\def\H{\mathcal H}
\def\N{\mathcal N}
\def\A{\mathcal A}
\def\one{\mathbf 1}

\def\eps{\epsilon}
\def\ignore#1{\relax}

\ignore{
\input BoxedEPS
\SetRokickiEPSFSpecial
\HideDisplacementBoxes
\SetEPSFDirectory{./graphics/} }

\input BoxedEPS
\SetTexturesEPSFSpecial
\HideDisplacementBoxes
\SetEPSFDirectory{:graphics:}

\begin{document}
%\frontmatter

\title[Crystal Bases and Kazhdan-Lusztig Polynomials]{Crystal Bases of
Quantum Affine Algebras and Affine Kazhdan-Lusztig Polynomials }

\author{Frederick M. Goodman and Hans Wenzl}

\address{Department of Mathemematics\\ University of Iowa\\ Iowa City, Iowa}

\email{goodman@math.uiowa.edu}

\address{Department of Mathematics\\ University of California\\ San Diego,
California}

\email{wenzl@brauer.ucsd.edu}

\bigskip

\begin{abstract} We  present a fast 
version of the algorithm of Lascoux, Leclerc, and Thibon for
the  lower global crystal base for the Fock representation of
quantum affine $\mathfrak{sl}_n$.  We also show that the coefficients
of the   lower global crystal base coincide with certain affine
Kazhdan-Lusztig polynomials.  Our algorithm allows fast computation
of decomposition numbers for tilting modules for quantum $\mathfrak{sl}_k$
at roots of unity, and for the Hecke algebra of type $A_n$ at roots of
unity.
\end{abstract}

\maketitle

\section{Introduction}

Lascoux, Leclerc, and Thibon in ~\cite{LLT} conjectured a connection between
the decomposition numbers for the Hecke algebra $H_n(q; \C)$ of type
$A_{n-1}$, for
 $q$  a primitive $l$-th root of unity in $\C$,  and the  {\em  lower global
crystal base}  for the simple highest weight module $M(\Lambda_0)$ of the
affine quantum group $\usl$. More explicitly, they conjectured that the
coefficients $\dlm(v)$ of the
 lower global crytal base, which {\em a priori} lie in $\Q[v]$, satisfy

\medskip
\noindent (1)\quad $\dlm(1) = \dlm$, where $\dlm = (S^\la: D^\mu)$ is the
multiplicity of the simple
$H_n(q; \C)$  module $D^\mu$ in a composition series for the Specht module
$S^\la$.

\medskip
\noindent(2) \quad $\dlm(v) \in \nat[v]$ (positivity).

\medskip\noindent The Specht module is a canonical indecomposable module of
$H_n(q; \C)$, which has, when $\la$ is $l$-regular, a unique simple quotient
$D^\la$; see ~\cite{Dipper-James1, Dipper-James2}. The first of these
conjectures was proved by Ariki ~\cite{Ariki}.
 Varagnolo and Vasserot ~\cite{Varagnolo-Vasserot} have recently proved a
more general conjecture of Leclerc and Thibon  ~\cite{LT} concerning
$q$-Schur algebras, and have  obtained a proof of the positivity conjecture
(2) as well.  A recent review of Geck ~\cite{Geck} discusses the work of LLT
and Ariki.
 
Lascoux, Leclerc, and Thibon also gave a recursive algorithm for computing
the  lower global crystal base in ~\cite{LLT}; we will consider 
improvements of their algorithm in this paper.

A similar connection between  certain affine Kazhdan-Lusztig polynomials and
decomposition numbers for the quantum group $U_q(\g)$, where $\g$ is any
simple Lie algebra over
$\C$ and $q$ is a root of unity, was  conjectured by Soergel in
~\cite{Soergel1} and  proved in ~\cite{Soergel2}.   Soergel's theorem is that
the affine Kazhdan-Lusztig polynomials 
${ n}_{\la, \mu}(v)$ described in ~\cite{Soergel1} satisfy
$$ { n}_{\la + \rho, \mu+\rho}(1) = (T_\mu: \nabla_\la),
$$ where $T(\mu)$ is the indecomposable
tilting module of $U_q(\g)$ with
highest weight
$\mu$ and $\nabla(\la)$ is the ``good" module with highest weight $\la$.
We quickly review the relevant definitions from [An1]:  A good module
$\nabla(\mu)$ is one induced from the character
$\mu$ of the Borel subalgebra $U_q^0 U_q^+$.  A Weyl module $\Delta(\la)$ is
$\nabla(-w_0 \la)^*$, where $w_0$ is the longest element of the Weyl group
of $\g$.  A tilting module is one with a filtration by  good modules and by
Weyl modules; there exists a unique indecomposable tilting module
$T(\mu)$ with highest weight $\mu$.  Finally
$(T_\mu:
\nabla_\la)$ is the  number of subquotients isomorphic to $\nabla(\la)$ in a
good filtration of $T(\mu)$.

The affine Kazhdan-Lusztig polynomials are also computed by a recursive
algorithm, which is described in  ~\cite{Soergel1}.

In this paper, we  present a fast version of the algorithm of ~\cite{LLT}
for the coefficients of the  lower global crystal base.  We also show that
the two polynomial analogues of decomposition numbers coincide (for Lie type
A), 
$$
\dlm(v) = { n}_{\la + \rho, \mu+\rho}(v). 
$$ 

The first idea of our modified LLT algorithm is to  fix an integer $k$ and
compute the $\dlm(v)$ only for Young diagrams $\la, \mu$ with no
more than $k$ rows; in this way, we can
identify Young diagrams with dominant integral weights of
$U_q(\mathfrak{sl}_k)$.    We wish to compute a formal sum:
$$
\tilde G(\mu) = \sum \{\dlm(v) \la  :  \ell(\la) \le k\}.
$$ 
The method of the original LLT algorthm is the following:  One considers
a certain path in the Weyl chamber of $\mathfrak{sl}_k$ from $\rho$ to
$\mu + \rho$, and all ``conjugate" paths obtained by
reflections 
in  hyperplanes of the affine Weyl group of $\mathfrak{sl}_k$
acting at level
$l$. (Here $\rho$ is one-half the sum of positive roots.)
The conjugate paths, which  end 
 at diagrams $\la + \rho$  which are in the orbit of
$\mu +
\rho$ under the action of the affine Weyl group, are   assigned
Laurent polynomial weights 
according to a rule involving the local geometry at points where the
path meets or leaves reflection hyperplanes.  
These weights provide a  first approximation to
$\tilde G(\mu)$:
$$
\tilde A(\mu) = \sum \{\alm(v) \la  :  \ell(\la) \le k\}.
$$ 
The next step is a  ``trimming" operation:  The approximation is
corrected by subtracting certain  multiples of $\tilde G(\la)$ for
$\la + \rho$ an endpoint of a conjugate path. 

In our modified algorithm, we take instead a path to  $\mu + \rho$
which begins at a nearby point  
 $\nu +
\rho$ which is as singular as possible, that is,  has large isotropy for the
action of the affine Weyl group. The recursive
computation  of $\tilde G(\mu)$ then involves 
$\tilde G(\nu)$ 
 as well as the  $\tilde
G(\la)$ for
$\la + \rho$ the endpoint of a conjugate path. However, as the path is
short, it has few conjugates; furthermore,  $\tilde G(\nu)$ is
relatively easy to compute because of the large isotropy group of
$\nu + \rho$.

When the point $\mu + \rho$ is sufficiently far from the walls
of the Weyl chamber, our algorithm (which was suggested by
our previous work [GW]) is roughly equivalent to Soergel's
algorithm [S1]
for affine Kazhdan-Lusztig polynomials (which involves
the periodic Kazhdan-Lusztig module). The main advantage of our
approach comes for diagrams which are near  the boundary of the 
Weyl chamber.

Regarding the equality of the polynomial analogues of the decomposition
numbers, note
that Schur-Weyl duality [Du]  already suggests that the
decomposition numbers coincide:
$
\dlm = (T_\mu: \nabla_\la)
$.   
Our proof of the equality 
$
\dlm(v) = { n}_{\la + \rho, \mu+\rho}(v)
$ 
is combinatorial;  we show that the $\dlm(v)$ satisfy the same
recursive relations which define the Kazhdan-Lusztig polynomials.
Varagnolo and Vasserot have recently given a
representation theoretic proof that the
$\dlm(v)$ coincide with Kazhdan-Lusztig
polynomials ~\cite{Varagnolo-Vasserot}, Theorem 9.3.
Steen Ryom-Hansen has informed us that he has also proved
the coincidence of the two polynomial analogues of decomposition numbers 
~\cite{Hansen}.

We intend to address generalizations of  this work
 to other Lie types, as well as relations to tensor ideals for
tilting modules, in future publications.

\section{Preliminaries} \label{section preliminaries}

Fix a positive integer $l$. Write $[[n]] = \Mod(n,l) = n - [n/l]l$; thus $n
= [n/l]l + [[n]]$, and
$0 \le [[n]] \le l-1$

We adopt the conventions of [M] regarding partitions,  Young diagrams and
tableaux. The Young diagram of a partition
$\la$ is the set of points $(i,j) \in \nat \times \nat$ such that
$1 \le j \le \la_i$; we identify the partition with its diagram.

For a point $(i,j) \in \nat \times \nat$, the {\em content} of $(i,j)$ is
$i-j$, and the {\em $l$-residue} is $[[i-j]]$; in particular these
quantities are defined for the nodes of a Young diagram. A Young diagram is 
called $l$-regular if it has no more than $l-1$ rows of any length. The {\em
length} $\ell(\la)$ of $\la$ is the number of (non-empty) rows. A  tableau
is regarded alternatively as a numbering of the nodes of a Young diagram, or
as an increasing sequence of Young diagrams.

A  node  $(i, \la_i) \in \la $ is  called {\em removable} if $i = \ell(\la)$
or if 
$i < \ell(\la)$ and 
$\la_{i+1} < \la_i$;  that is,  if one obtains a Young diagram by removing
the node from $\la$.  The node is called a {\em removable $r$-node} if it is
removable and if its $l$-residue is equal to $r$.  

A  point $(i, \la_i +1)$ is called an {\em indent} node of $\la$ if $i = 1$ 
or if $i > 1$ and
$\la_{i-1} > \la_i$, that is, if  one obtains a Young diagram by adding a
new node to $\la$ at $(i, \la_i +1)$.   Note that $(\ell(\la) + 1, 1)$
counts as an indent node. The node is called an {\em indent $r$-node} of
$\la$ if it is  an indent node and if its $l$-residue is equal to $r$.

We will frequently use the notion of a  {\em face} of a closed convex set of
a Euclidean space $\R^k$.  A face $F$  of a closed convex set $C$ is an
extreme convex subset of $C$; this means that if a point $f \in F$ is a
convex combination of points $a, b \in C$, then $a, b \in F$.   The {\em
boundary} $\partial F$ of a face $F$ is the union of all proper subfaces;
this does not generally coincide with the topological boundary of $F$ in
$\R^k$. The {\em interior} of a face
$F$ is $F \setminus \partial F$.  An open face is the interior of a closed
face.

We will denote the unit vector in $\R^k$ with a 1 in the $i\th$ coordinate
and zeros elsewhere by $\eps_i$.

\section{The LLT algorithm}

We first recall the setting for the algorithm of Lascoux, Leclerc, and
Thibon. The affine quantum universal enveloping algebra $\usl$, defined over
$\Q(v)$, has generators $e_i$ and $f_i$, for $0 \le i \le l-1$, and $v^h$,
for
$h$ in the dual weight lattice $P\,\check{}$ of
$\widehat{\mathfrak{sl}}_l$.  We refer, for example to ~\cite{LLT}, Section
4, for the list of relations satisfied by the generators.  Write
$U^-$ for the subalgebra generated by  the $f_i$.

The algebra
$\usl$ has a representation on the ``Fock space" $\mathcal F$, which is 
the $\Q(v)$-vector space spanned by Young diagrams of all sizes.  The
generators
$q^h$ are diagonal on the basis of Young diagrams.  The generators $e_i$ act
as ``annihilation operators,"  removing nodes from a Young diagram, and the
generators
$f_i$ act as ``creation operators," adding nodes to a Young diagram. This
representation was described by Hayashi [H] and Misra and Miwa [MM]; see also
~\cite{LLT}, Section 4.  We will need  the explicit description of the
action only for the $f_i$:

\begin{equation} \label{equation formula for fi} f_i \la = \sum_\nu
v^{N(\la, \nu)} \nu,
\end{equation} where the sum is over all diagrams $\nu$ which can be
obtained from $\la$ by adding one node of $l$-residue i, and $N(\la, \nu)$
is given as follows. If $\nu$ is obtained from $\la$ by adding a node of
$l$-residue $i$ in a certain row $r$, then $N(\la, \nu)$ is the number of
indent $i$-nodes of $\la$ in rows
$r' < r$, less the number of removable  $i$-nodes of $\la$ in rows
$r' < r$.  Of course, $f_i \la = 0$ if $\la$ has no indent $i$-nodes.

We write $M$ for the cyclic submodule of $\mathcal F$ generated by the empty
diagram $\emptyset$, $M = \usl \emptyset = U^-  \emptyset$. Then $M$ is
isomorphic to the simple integrable highest weight module
$M(\Lambda_0)$.  Let $A$ be the subring of $\Q(v)$ of rational functions
with no poles at $v = 0$, and let $L$ denote the $A$-span of all Young
diagrams. Let $\umq$ denote the $Q[v,v\inv]$-subalgebra of $U^-$ generated by
all divided powers 
$f_i^{(m)} = f_i^m/[m]!$. Here $[m]$ is the $v$-integer,
$[m] = \displaystyle
\frac{v^m - v^{-m}}{v - v\inv}$, and $[m]! = [m] [m-1] \cdots [1]$. Let
$M_\Q$ denote the cyclic $\umq$ module, $M_\Q = \umq \emptyset$.

There is an involution $a \mapsto \bar a$ of $\usl$ such that the $e_i$ and
$f_i$ are self-dual, $\bar v = v\inv$, and $(v^h)^- = v^{-h}$. This induces
an involution on $M$ by $(a \emptyset)^- =  \bar a \emptyset$.

The following theorem is due to Kashiwara (see ~\cite{Kashiwara1,
Kashiwara2}).

\medskip
\begin{theorem}  There is a unique $\Q[v, v\inv]$-basis $\{G(\mu)\}$  of
$M_\Q$
 indexed by $l$-regular Young diagrams $\mu$, and satisfying:

\medskip\noindent{\rm (G1)}\quad  $G(\mu) \equiv \mu \quad (\mod\ L)$; and

\medskip\noindent{\rm (G2)}\quad $G(\mu) $ is self-dual. 
\end{theorem}

\medskip This basis is called the {\em lower global crystal basis}.  The
algorithm  given by Lascoux, Leclerc and Thibon computes the elements of
this basis.

As observed in ~\cite{LLT}, in order to obtain the lower global crystal
basis of
$M(\Lambda_0)$, it is convenient to first compute for each
$l$-regular Young diagram $\mu$ an element $A(\mu) \in M_\Q$ satisfying:

\medskip
\noindent (A1)\quad $A(\mu) = \sum_\la \alpha_{\la \mu} \mu$, where
$\alpha_{\mu \mu} = 1$, 
$\alpha_{\la \mu} \in \Z[v,v\inv]$, and $\alpha_{\la \mu} = 0$ unless
$\la \trianglelefteq \mu$, and

\medskip
\noindent (A2)\quad $A(\mu)$ is self-dual. 
\medskip

\noindent Such elements  $(A_\mu)$ are by no means unique.  

The $G(\mu)$  can be  computed recursively from the $A(\mu)$ by a triangular
reduction algorithm (Gaussian elimination). First, if $\mu$ is the smallest
$l$-regular Young diagram of a given size $n$ in lexicographic order, then
$G(\mu) = A(\mu) = \mu$, as follows from (A1)  and (G1). Next, given an
$l$-regular Young diagram of size $n$, one can compute
$G(\mu)$ using $A(\mu)$ and the $G(\mu')$ for $\mu' < \mu$ in lexicographic
order.  Namely, if there is some $\mu' < \mu$  such that the coefficient
$\alpha_{\mu' \mu}(v) \not\in v\Z[v]$, then take the lexicographically
largest such $\mu'$.  Necessarily $\mu'$ is $l$-regular.  Writing
$\alpha_{\mu' \mu}(v) = \displaystyle\sum_{n = -N}^N \alpha_{\mu' \mu}^n
v^n$, put
$\gamma_{\mu'} = \displaystyle \sum_{n = -N}^{-1} \alpha_{\mu' \mu}^n (v^n +
v^{-n}) + 
\alpha_{\mu' \mu}^0$, and replace $A(\mu)$ by $A(\mu) - \gamma_{\mu'}
G(\mu')$. Continue in this way until all coefficients in $A(\mu)$ for
diagrams $\mu' < \mu$ lie in $v\Z[v]$.  Then one has $A(\mu) = G(\mu)$ by
uniqueness of $G(\mu)$.

Next we consider how to compute elements $A(\mu)$ satisfying (A1) and (A2).

\medskip
\begin{definition} Let $T$ be a standard (skew) tableau of size $n$.  The
{\em
$l$-residue sequence} associated to $T$ is described as follows. For each
$i$, $1 \le i \le n$, let $a_i$  be the $l$-residue of  the node
containing $i$ in the tableau $T$. Now write
\begin{equation} (a_1, \dots, a_n) = (r_1^{m_1}, r_2^{m_2},\dots, r_s^{m_s}).
\end{equation} More precisely,  consider the maximal constant subsequences
of the sequence
$(a_i)$ and  let $r_i$ be the constant value of the $i\th$ constant
subsequence and $m_i$ be the length of the $i\th$ constant subsequence.
\end{definition}

\medskip
\begin{definition}  We say that two standard tableaux are  {\em
$l$-conjugate} if they have the same $l$-residue sequence.  Similarly we say
that two standard skew tableaux are $l$-conjugate if they have the same
starting shape and the same $l$-residue sequence.
\end{definition}

\medskip Note that $l$-conjugate (skew) tableaux can be of different shapes. 

\medskip

\begin{definition} \label{foft} Let $T$ be a standard (skew) tableau.  
Define a corresponding element $f(T) \in U_\Q^-$ by
\begin{equation} f(T) = f_{r_s}^{(m_s)}\cdots f_{r_1}^{(m_1)},
\end{equation} where $(r_1^{m_1}, r_2^{m_2},\dots, r_s^{m_s})$ is the
$l$-residue sequence of  $T$.  
\end{definition}

\medskip
\begin{definition} If $T$ is a standard tableau (i.e. not a skew tableau),
define $A(T) = f(T) \emptyset$.
\end{definition}

\medskip Note that the element $A(T)$ is always  self-dual.

\medskip The following combinatorial lemma from ~\cite{LLT} is important for
understanding the properties of the elements of the form $A(T)$.  

\medskip
\begin{lemma} \label{lemma action of divided power} For any Young diagram
$\la$,
$$ f_r^{(m)} \la = \sum_\mu v^{N(\la,\mu)} \mu,
$$ where the sum is over all Young diagrams $\mu$ such that $\mu \setminus
\la$ consists of $m$ nodes all with $l$-residue $r$, and 
$N(\la, \mu)$ is computed as follows:  Let 
$$ T_0 = (\la = \mu^{(0)} \subseteq \mu^{(1)} \subseteq \cdots \subseteq
\mu^{(m)} = \mu)
$$ be the standard skew tableau of shape $\mu\setminus\la$ in which the $m$
nodes of $\mu\setminus\la$ are filled from top to bottom. Put
$$ N(T_0) = \sum _i N(\mu^{(i-1)},\mu^{(i)}),
$$ where now $N(\mu^{(i-1)},\mu^{(i)})$ is as in the explanation following
Equation \ref{equation formula for fi}.  Then
$$ N(\la,\mu) = N(T_0) + \binom{m}{2}.
$$
\end{lemma}

\begin{proof}  See the proof of Lemma 6.2 in ~\cite{LLT}.
\end{proof}

\medskip It is evident that $N(\la, \mu)$ depends only on the position of the
$m$ nodes of $\mu \setminus \la$ relative to other possible indent $r$-nodes
of $\la$ and possible removable $r$-nodes of $\la$.

\medskip
\begin{lemma} \label{lemma property L tabs1}
 Suppose that every removable $r$-node of
$\la$ lies below every indent $r$-node of $\la$ and furthermore that $\mu
\setminus
\la$  consists of the first (highest) $m$  indent $r$-nodes of $\la$. Then
$N(\la, \mu) = 0$.
\end{lemma} 

\begin{proof}  In this case $N(T_0) = - \binom{m}{2}$; compare the proof of
Lemma 6.4 in ~\cite{LLT}.
\end{proof}

\medskip We describe a map from standard (skew) tableaux to certain (row-
and  column-strict) semi-standard (skew) tableaux of the same shape,  as
follows.   Let $T$ be a standard (skew)  tableau of shape $\la\setminus\nu$,
with $l$-residue sequence $(r_1^{m_1}, r_2^{m_2},\dots, r_s^{m_s})$.  The
corresponding semi-standard tableau is
$$ S(T) = (\nu = \la^{(0)} \subseteq \la^{(1)} \subseteq \cdots \subseteq
\la^{(s)} =
\la),
$$ where $\la^{(i)}\setminus\nu$ is the subdiagram of $\la\setminus \nu$
containing the numbers 1 through $m_1 + \cdots + m_i$.  Thus $\la^{(i)}
\setminus \la^{(i-1)}$ is a skew diagram consisting of $m_i$ nodes all with
$l$-residue equal to
$r_i$.  As there are $m_i\,!$ different standard fillings of
$\la^{(i)} \setminus \la^{(i-1)}$ for each $i$, 
$S(T)$ is the image of $\prod_i m_i!$ different standard (skew) tableaux. 
Define
\begin{equation}
\wt(T) =  \wt(S(T)) = v^{\sum N(\la^{(i-1)}, \la^{i)})},
\end{equation} where $N(\la^{(i-1)}, \la^{i)})$ is as in Lemma 
\ref{lemma action of divided power}

The operator $f(T)$ applied to the empty diagram may be thought of as
generating all semi-standard tableaux $S(T')$, where $T'$ is a standard
tableau $l$-conjugate to $T$. Then 
$$A(T) = \sum_S \wt(S) \la(S),$$ where  the sum is over all $S = S(T')$ with
$T'$ $l$-conjugate to $T$, and
$\la(S)$ denotes the shape of the semi-standard tableau $S$.

Let $T$ be a standard (skew) tableau of shape $\la\setminus\nu$,  let
$(r_1^{m_1}, r_2^{m_2},\dots, r_s^{m_s}) $ be the $l$-residue sequence of
$T$ and let
$$ S(T) = (\nu = \la^{(0)} \subseteq \la^{(1)} \subseteq \cdots \subseteq
\la^{(s)} =
\la)
$$ the semi-standard tableau associated to $T$. 
We require a technical condition on the tableau $T$ which will insure
that $A(T)$ will satisfy properties {\rm (A1)} and {\rm (A2)}.
The following is a sufficient condition:

\medskip\noindent (L) 
\quad   For every $i$, each removable
$r_i$ node of $\la^{(i-1)}$ lies below each indent $r_i$ node of
$\la^{(i-1)}$, and $\la^{(i)} \setminus \la^{(i-1)}$ consists of the first
$m_i$ indent
$r_i$-nodes of $\la^{(i-1)}$.

\medskip
\vbox{
\begin{lemma} \label{lemma property L tabs2} Let $T$ be a standard tableau
with property {\rm (L)} and with shape
$\mu$.
\begin{enumerate}
\item  If $T'$ is  a standard tableau $l$-conjugate to $T$, with shape
$\mu'$, then $\mu' \trianglelefteq \mu$.
\item  $A(T)$ satisfies properties {\rm (A1)} and {\rm (A2)}.
\end{enumerate}
\end{lemma} }

\begin{proof}  Point (a) is immediate from the definition of property (L).
All elements of the form $A(T)$ satisfy (A2).  Finally property (A1) follows
from Lemmas \ref{lemma action of divided power} and  
\ref{lemma property L tabs1}.
\end{proof}

\medskip Lascoux, Leclerc, and Thibon describe for each $l$-regular Young
diagram $\mu$ a certain standard tableau $T(\mu)$ of shape $\mu$ which has 
property (L), namely the standard tableau which fills $\mu$ from top to
bottom along ``$l$-ladders."  So putting 
$$ A(\mu) =  A(T(\mu)) 
$$ for each $l$-regular Young diagram $\mu$ gives a family $(A(\mu))_\mu$
with the desired properties, by Lemma \ref{lemma property L tabs2}.

\section{A Modified Algorithm}

In this section, we explain a fast variant  of the LLT algorithm. The LLT
algorithm is relatively slow because the  $A(\mu)$ produced in the first
step tend to have contributions from many diagrams
$\mu'$ which do not contribute to $G(\mu)$; this means that the recursive
computation during the triangular reduction step becomes very complicated.
So the key to a more efficient computation is to find a way to  compute
$A(\mu)$ which have fewer extraneous contributions, and so to reduce the
complexity of the recursion.

We fix an integer $k$, and consider only Young diagrams of length no more
than $k$.  So for an $l$-regular Young diagram $\mu$, we consider the problem
of computing the truncation of $G(\mu)$
$$
\tilde G(\mu) = \sum \{\dlm(v) \la  :  \ell(\la) \le k\}.
$$ It is more or less clear that the LLT algorithm, modified only by 
throwing out diagrams with more than $k$ rows, computes the $\tilde G(\mu)$.
Formally, we note that the span $\mathcal F_k'$ of diagrams with more than
$k$ rows is a  $\um$-submodule of $\mathcal F$.  So we work in  the quotient
module $\mathcal F_k = \mathcal F/\mathcal F_k'$, and the cyclic submodule 
$M_k = M/(M \cap \mathcal F_k') $. The quotient module $\mathcal F_k$ has
its natural basis labelled by Young diagrams with no more than $k$ rows; and
we regard the 
$\tilde G(\mu)$ as elements of $M_k$. Note that $M \cap \mathcal F_k'$ is
also invariant under the involution of $M$, because the $G(\mu)$ with
$\ell(\mu) > k$ is a self-dual basis, so
 the involution also passes to the quotient
$M_k$.  Let $L_k \subseteq \mathcal F_k$ denote the $A$-span of diagrams
with no more than $k$ rows, i.e. the image of the $A$-lattice $L$ in 
$\mathcal F_k$, and let $M_{\Q, k}$ denote the image of $M_\Q$ in 
$\mathcal F_k$.

\medskip
\begin{lemma}  The set of $\tilde G(\mu)$ where $\mu$ is an $l$-regular
Young diagram with no more than $k$ rows is the unique 
$\Q[v, v\inv]$  basis of
$M_{\Q,k}$ satisfying:

\medskip\noindent{\rm ($\tilde {\rm G}$1)}  $\tilde G(\mu) \equiv \mu \quad
(\mod\  L_k)$

\medskip\noindent{\rm ($\tilde {\rm G}$2)}  $\tilde G(\mu)$ is self-dual.
\end{lemma}

\begin{proof}  The set of $\tilde G(\mu)$ does have the stated properties.
For uniqueness, if $\{\tilde H(\mu)\}$ is another such basis, then each
$\tilde H(\mu)$ has a self-dual pre-image $ H(\mu)$ in $M_\Q$, and the set
of $H(\mu)$ together with the set of $G(\mu)$ with $\ell(\mu) > k$ is then a
basis of $M_\Q$ satisying (G1) and (G2).  By uniqueness of the lower global
crystal basis, one has $H(\mu) = G(\mu)$, and hence
$\tilde H(\mu) = \tilde G(\mu)$.
\end{proof}

\medskip In order to compute the $\tilde G(\mu)$, it will suffice to find
elements
$\tilde A(\mu) \in M_k$, for $l$-regular diagrams of length no more than $k$,
satisfying

\medskip
\noindent ($\tilde {\rm A}$1)\quad $\tilde A(\mu) = \sum\{ \alpha_{\la
\mu}(v) \la : \ell(\la)
\le k\}$, where
$\alpha_{\mu
\mu} = 1$, 
$\alpha_{\la \mu} \in \Z[v,v\inv]$, and $\alpha_{\la \mu} = 0$ unless
$\la \trianglelefteq \mu$; and

\medskip
\noindent ($\tilde {\rm A}$2)\quad $\tilde A(\mu)$ is self-dual. 
\medskip

The $\tilde G(\mu)$ can then be obtained from the $\tilde A(\mu)$ by
triangular reduction, as before.

\medskip We now consider how to compute such elements $\tilde A(\mu)$.
Recall that the
$l$-residue sequence, and the operator
$f(T)$ are defined for standard skew tableaux as for standard tableaux.  
Henceforth, all Young diagrams will have  no more than $k$ rows, so we will
no longer mention this.

\medskip
\begin{definition}  If $T$ is a standard skew tableau of some shape $\mu
\setminus \nu$,  define
$
\tilde A(T) = f(T) \tilde G(\nu).
$
\end{definition}
\medskip The proof of the following lemma is left to the reader.
\medskip
\begin{lemma}\label{lemma property L tabs3} Let $T$ be a standard skew 
tableau with property {\rm (L)} and with shape
$\mu\setminus\nu$.
 Then  $\tilde A(T)$ satisfies properties  {\rm ($\tilde {\rm A}$1)} and 
{\rm ($\tilde {\rm A}$2)}.
\end{lemma}

\medskip We identify Young diagrams with dominant   integral weights of 
$\mathfrak{gl}_k$, by 
$$
\la \mapsto \sum_{i = 1}^{k-1} (\la_i - \la_{i+1}) \Lambda_i + \la_k
\Lambda_k,
$$  where the $\Lambda_i$ are the fundamental weights. The fundamental
weight  
 ($1 \le i \le k$),  corresponds to the diagram with $i$ nodes arranged in
one column.  The half-sum $\rho$ of positive roots of $\mathfrak{sl}_k$
corresponds to the partition
$(k-1, k-2, \dots, 1, 0)$.   The positive Weyl chamber $\mathcal C$ of
$\mathfrak{sl}_k$ is the cone generated by $\{\Lambda_i: 1 \le i \le k-1 \}$
in the dual $\mathfrak h^*$ of the diagonal subalgebra.

For an element $A = \sum_\la \alpha_\la \la \in \mathcal F_k$, let us write
$\Lambda_k + A = \sum_\la \alpha_\la \ (\la + \Lambda_k)$.

\begin{proposition} \label{proposition adding Lambda_k} If both $\mu$ and
$\mu +
\Lambda_k$ are
$l$-regular, one has
$$\tilde G(\mu + \Lambda_k) = \Lambda_k + \tilde G(\mu).$$
\end{proposition}

\begin{proof}  We prove this lemma here only under the assumption that
$k < l$.  In this case, all Young diagrams are $l$-regular. The general case
follows from Proposition 
\ref{proposition independence of l}.

It follows from the LLT algorithm for the lower global crystal base that the
diagrams $\la$ appearing in the expansion of $\tilde G(\mu)$ have the same
$l$-core as $\mu$ and satisfy $\la \trianglelefteq \mu$. Therefore if $\mu$
is the lexicographically smallest diagram of a given size with a given $l$
core, then $\tilde G(\mu) = \mu$.  The same holds for
$\mu + \Lambda_k$, so $\tilde G(\mu + \Lambda_k) = \mu + \Lambda_k =
\Lambda_k + \tilde G(\mu)$.   This argument applies in particular to  the
empty diagram.

Now let $\mu$ be a non-empty diagram and suppose inductively that the
assertion holds for all diagrams with fewer nodes, and for all diagrams with
the same number of nodes which are lexicographically smaller. Let $\mu'$ be
the diagram obtained by removing the highest removable node of $\mu$  (i.e.,
that with the least row index), let $r$ be the residue and $i$ the row index
of that node.  By the induction, hypothesis, we have
$
\tilde G(\mu' + \Lambda_k) = \Lambda_k + \tilde G(\mu')
$. The diagram $\mu'$ has no indent $r$-nodes and no removable $r$-nodes in
rows $j < i$, so we can take $\tilde A(\mu) = f_r \tilde G(\mu')$. For the
same reason, 
$$\tilde A(\mu + \Lambda_k ) = f_{r+1} \tilde G(\mu' + \Lambda_k) = f_{r+1}
( \Lambda_k +\tilde G(\mu') ) = \Lambda_k + f_r \tilde G(\mu)=
\Lambda_k + \tilde A(\mu).
$$ The conclusion now follows from the induction hypothesis.
\end{proof}

\medskip
\begin{corollary} If $\mu$ is an $l$-regular Young diagram with $\mu_k \ne
0$, then
$$\tilde G(\mu) = \mu_k\Lambda_k + \tilde G(\mu - \mu_k \Lambda_k).$$ 
\end{corollary}

It now suffices to give an algorithm for computing $\tilde G(\mu)$ when
$\mu$ is an $l$-regular diagram with $\mu_k = 0$. Our strategy will be to
choose for any $l$-regular Young diagram $\mu$ an $l$-regular Young diagram
$\nu \subseteq \mu$ and a standard skew tableau 
$T$ of shape
$\mu \setminus \nu$, and to define $\tilde A(\mu) = \tilde A(T) = f(T) \tilde
G(\nu)$.  In the original  LLT algorithm, one always takes $\nu$ to be the
empty diagram and $T$ to be the standard ladder tableau.  Here we will choose
$\nu$ and $T$ in a way which makes the recursive computation of the $\tilde
G(\mu)$ more efficient.

\medskip
\noindent{\bf Case 1. Critical and Interior Diagrams.}
\medskip

A Young  diagram $\mu$  is called  {\em $k$-critical} if $\mu_i-\mu_{i+1} +
1 $  is divisible by $l$ for $i=1,2,\ ...\ k-1$. Equivalently, $(\mu + \rho,
\alpha)$ is divisible by $l$ for all roots $\alpha$ of $\mathfrak{sl}_k$. 
The Steinberg weight $(l-1)\rho$ is the smallest
$k$-critical dominant integral weight.  We call a Young diagram $\mu$ {\em
interior} if it lies above the Steinberg weight, i.e.,  $\mu  - (l-1)\rho$
is a dominant integral weight. All interior diagrams are $l$-regular, and
all $k$-critical diagrams are interior. 

The {\em fundamental box} $B$ is the set
$$
\{ \sum_{i = 1}^{k-1}  m_i \Lambda_i : 0 \le m_i < l \}.
$$ The set of interior Young diagrams
 is tiled by translates of the fundamental box $B$ by $k$-critical Young
diagrams; i.e. for each interior Young diagram $\mu$, there is a unique
$k$-critical Young diagram $\mu_c$ such that
$\mu \in \mu_c + B$.

\medskip
\begin{proposition} \label{Gmu=mu} If $\mu$ is a $k$-critical diagram,  then
$\tilde G(\mu) =
\mu$.
\end{proposition}

\begin{proof}  It is possible to give a simple combinatorial proof, adapting
arguments from ~\cite{Goodman-Wenzl}.  Here we will make use of Ariki's
theorem, and the positivity result of Varagnolo and Vasserot, in the
interest of brevity.  

\medskip It is well known that the decomposition numbers for the Hecke
algebra $H_n(q; \C)$ satisfy
$\dlm = 0$ for all Young diagrams $\la \ne \mu$ with no more than $k$-rows;
see, for example ~\cite{Goodman-Wenzl}.  It follows from $\dlm(v) \in
\nat[v]$ and
$\dlm(1) = \dlm$, that $\dlm(v) = 0$ as well  for all Young diagrams $\la
\ne \mu$ with no more than $k$-rows.
\end{proof}

%%%%%%%%%%%%%%%%%%%%%%%%%%%%%%%%%%%%%%%%%%%%%%%%%%%%%%%%%%%%%%%%%%%
%%%%%  THE FOLLOWING IS AN ELEMENTARY PROOF WHICH IS NO LONGER NEEDED.
%%%%%%%%%%%%%%%%%%%%%%%%%%%%%%%%%%%%%%%%%%%%%%%%%%%%%%%%%%%%%%%%%%%

\ignore{Consider first the lexicographically lowest
 $k$-critical diagram $\mu$ of any given size $n$. Let $t$ be the standard
filling of $\mu$ from top to bottom along $l$-ladders.  By ~\cite{LLT},
Lemma      , $A(t) = \mu + \sum_\la \alm(z) \la$, where the sum goes over
$\la$ which satisfy $\la < \mu$, and which are the endpoints of tableaux
equivalent to $t$. But the endpoints of tableaux equivalent to
$t$ must also be $k$-critical diagrams of size $n$; since $\mu$ is the
lowest  such diagram, it follows that  
$A(t) = \mu$.  Hence  $G(\mu) = \mu$.  This conclusion holds in particular
for the smallest $k$-critical diagram 
$\mu^0 = [(k-1)l, (k-2)l, \dots, l, 0]$. }

\ignore{ Next consider a $k$-critical diagram $\mu \ne \mu^0$ which is {\em
reduced}, that is $\mu_k = 0$.  Assume that $G(\mu') = \mu'$ for all
lexicographically lower
$k$-critical diagrams $\mu'$  of the same size and all $k$-critical diagrams
of smaller size.  Note that $\mu + \rho$ has the form
$[a_1 l, a_2 l, \dots, a_{k-1}l, 0]$, where $a_1 > a_2 > \cdots > a_{k-1} >
0$. Let $j$ be the length of the Young diagram (paritition)  and $[a_1 -k+1
, a_2 -k + 2, \dots, a_{k-1}-1, 0] $. One can check that 
$j$ is also the number $i$ such that $\mu_i > \mu^0_i$. Furthermore $\mu' =
\mu - l e^{(j)}$ is also a  $k$-critical diagram. By the induction
assumption $G(\mu') = \mu'$. Now consider the element of $\ml0$  
\begin{equation} A =  f_0^{(j)}\dots f_2^{(j)} f_1^{(j)} G(\mu') = 
f_0^{(j)}\dots f_2^{(j)} f_1^{(j)} \mu'.
\end{equation} }

\ignore{ I claim that $A = \mu + \sum_\la \alm(z) \la$, and that if $\alm(z)
\ne 0$, then $\la < \mu$ and $\alm(z) = 1$.  Since $\mu'$ is $k$-critical,
it has an indent of $l$-residue 1 in each row.  Therefore $f_1^{(j)} \mu' =
\sum_e z^{\phi(e)} (\mu' + e)$,    where the sum ranges over all $k$-vectors
$e$ having $j$ entries equal to 1 and $k-j$ entries equal to 0.  Referring
to the computation in ~\cite{LLT}, Lemma     , one sees that the exponent
$\phi(e)$ is equal to the number of pairs $r < s$ such that $e_r = 0$ and
$e_s = 1$. Next, a diagram $\mu' + e$ has indents of $l$-residue 2 exactly
in the rows $r$ such that $e_r = 1$, so $f_2^{(j)} f_1^{(j)} \mu' = 
\sum_e z^{\phi(e)} (\mu' + 2e)$.  Continuing, one obtains that
$f_{l-1}^{(j)}\cdots f_2^{(j)} f_1^{(j)} \mu' = 
\sum_e z^{\phi(e)} (\mu' + (l-1)e)$.  Next one wants to add $j$ cells with
$l$-residue $0$.  If it is not possible to add these cells to $\mu' + (l-1)
e$, then $f_0^{(j)} (\mu' + (l-1) e) = 0$.  If it is possible to add these
cells to $\mu' + (l-1) e$,  then again  only in those rows $r$ such that
$e_r = 1$; but it this case the diagram
$\mu' + (l-1) e$ also has removable corners of $l$-residue $0$ in all rows
$r$ such that $e_r = 0$.  It follows that 
$f_0^{(j)} (\mu' + (l-1) e) = z^{-\phi(e)} (\mu' + l e)$.  So finally,
$f_0^{(j)}\dots f_2^{(j)} f_1^{(j)} \mu' = \sum_e (\mu' + l e)$ where the
sum is over all $e$ such that $\mu' + l e$ is a Young diagram. All such
diagrams are lexicographically less or equal to $\mu = \mu' + l e^{(j)}$.
This proves the claim.  }

\ignore{ Now using the induction hypothesis again, we have
$A = \mu + \sum_\la \la = \mu + \sum_\la G(\la)$, where the sum is over
certain diagrams $\la < \mu$.  Therefore  triangular reduction of $A$ as in
~\cite{LLT} gives $G(\mu) = \mu$. }

\ignore{ Next consider the case that $\mu$ is not reduced; say $\mu_k = s >
0$. Then $\mu' = \mu - s e^{(k)}$ is a reduced $k$-critical diagram, which
by  induction can be assumed to satisfy $G(\mu') = \mu'$.  If $[[\mu'_1]] =
r$, then  $f_{[[r+1]]}^{(s)} G(\mu') = \mu$.  Hence $G(\mu) = \mu$. }
%%%%%%%%%%   END OF DISCARDED PROOF  %%%%%%%%%%%%%%%%%%%%%%%%%%%%%%%%
%%%%%%%%%%%%%%%%%%%%%%%%%%%%%%%%%%%%%%%%%%%%%%%%%%%%%%%%%%%%%%%%%%%%%

Let $\mu$ be an interior Young diagram and let $\mu_c$ be the associated
critical diagram, such that $\mu \in \mu_c + B$.  Let
$\mu - \mu_c = \sum_{i = 1} ^{k-1} d_i \Lambda_i$, where $0 \le d_i < l$ for
each $i$. Consider the skew tableau $T$ from $\mu_c$ to $\mu$

\begin{equation}\begin{split} &\mu^{(k)} = \mu_c \rightarrow \cr
 &\mu^{(k)} + 
 \Lambda_{k-1} \rightarrow \dots \rightarrow
\mu^{(k)} +  d_{k-1} \Lambda_{k-1} = \mu^{(k-1)}  \rightarrow \cr
 &\mu^{(k-1)} + 
 \Lambda_{k-2} \rightarrow \dots \rightarrow
\mu^{(k-1)} +  d_{k-2} \Lambda_{k-2} = \mu^{(k-2)}  \rightarrow \cr &\vdots
\cr &	\mu^{(2)} + \Lambda_{1} \rightarrow \dots 
\rightarrow \mu^{(2)} + d_1 \Lambda_{1} = \mu,
\end{split}\end{equation} where at each stage the $i$ cells of $\Lambda_i$
are added from top to bottom.

It is easy to see that the skew tableau $T$ has property (L), so by Lemma 
\ref{lemma property L tabs2}
$\tilde A(\mu) = \tilde A(T) = f(T) \mu$ satisfies  ($\tilde{\rm  A}$1) and
($\tilde{ \rm A}$2).

\medskip
\begin{remark}  This construction was suggested by our approach to the
decomposition numbers  in ~\cite{Goodman-Wenzl}. For $\mu$ an interior
diagram, the definition of the element
$\tilde A(\mu)$  is a `$q$ -version' of the construction in section 5 of
~\cite{Goodman-Wenzl}.  There we defined $N(\la,\mu)$ to be the number of
paths from $\mu_c$ to $\la$ which are conjugate to the skew tableau $T$, and 
\begin{equation} n(\la, \mu) = \frac{N(\la, \mu)}  {N(\mu, \mu)} =
\frac{N(\la, \mu)}  {\prod_{i=1}^{k-1} (i!)^{d_i}}.
\end{equation} One has $n(\la, \mu) = \alm(1)$, where
$\tilde A(\mu) = \sum_\la \alpha_{\la \mu}(v) \la$.
  We showed in ~\cite{Goodman-Wenzl} that
$n(\la,
\mu)$ is an upper bound for the decomposition number $\dlm$.
\end{remark}

\medskip  Using section 5 of this paper, Soergel's theorem [S2], and results
on tensor ideals for tilting modules [O], one can show that the  algorithm
for $\tilde G(\mu)$, or a minor modification of it, for interior diagrams
$\mu$ has the properties:

\begin{enumerate}
\item  The only diagrams $\mu'$ such that $\tilde G(\mu')$ is used in the
recursive calculation of $\tilde G(\mu)$ are interior diagrams.

\item  The complexity of the algorithm for interior diagrams $\mu$ is
uniformly bounded, for fixed $k$ and $l$, independent of the size of the
diagram $\mu$. The dependence on $l$ can also be eliminated using
Proposition \ref{proposition independence of l}.
\end{enumerate}

We do not have a direct elementary proof of these statements.

\medskip
\noindent{\bf Case 2. Non-interior diagrams. }
\medskip

We now consider how to define $\tilde A(\mu)$ when $\mu$ is not an interior
diagram.  For simplicity, we assume that $k \le l$, so that all Young
diagrams are $l$-regular (aside from the diagrams $n \Lambda_k$, if $k = l$.)
  This assumption is not essential, and we will indicate afterwards how to
modify the procedure for $k > l$.

If $\mu$ is a non-interior Young diagram, then $\mu + \rho \in a + B$ for
some diagram $a$ located on one or more boundary hyperplanes of the positive
Weyl chamber $\mathcal C$, and satisfying $(a, \alpha_i)$ is divisible by
$l$ for all simple roots $\alpha_i$. 

If $a = 0$, i.e., if $\mu + \rho$ is contained in the fundamental box $B$,
then we compute $\tilde A(\mu)$ by the LLT algorithm; that is, take $T$ to
be the standard ladder tableau of shape $\mu$ and put
$\tilde A(\mu) = f(T) \emptyset$.

If $a \ne 0$, we proceed as follows:  Let $A^+$ denote the lowest $l$-alcove
in the positive Weyl chamber $\mathcal C$, namely
$$A^+ = \{ x :  x_1 > x_2 \ge \cdots >  x_k, \text{ and  }  x_1 - x_k  < 
l\}.$$ There is a unique closed face $F$ of $(a + A^+)^-$ of smallest
dimension such that $a \in F$, and the interior $F^0$
 of the face lies in the interior of the positive Weyl chamber.  In fact, let
$I = \{ i_1, \dots, i_s\} \subseteq \{1, 2, \dots, k-1\}$ be the complete
list of  indices $i$ such that the simple root 
$\alpha_i$ satisfies
 $(a, \alpha_i) = 0$.  Then $F$ is the convex hull of 
$\{a\} \cup \{a + l \Lambda_i : i \in I\}$.  Because we are assuming
$k \le l$,  there exist integer points on $F^0$; to be definite, take the
point $p = a + \Lambda_{i_1} +  \Lambda_{i_2} + \cdots + \Lambda_{i_s}$. We
have $\mu + \rho =  a+ \sum_i d_i \Lambda_i$, where $0 \le d_i < l$, and
necessarily $d_i >0$ for $i \in I$. Hence $\mu + \rho = p + \sum_i d_i'
\Lambda_i$, where $d_i' = d_i$ if
$i \not\in I$ and $d_i' = d_i -1$ if $i \in I$. Since $p$ is in the interior
of the Weyl chamber, it has the form $\nu + \rho$ for some Young diagram
$\nu$, and we have
$\mu = \nu + \sum_i d_i' \Lambda_i$.  Now we take the skew tableau of shape
$\mu \setminus \nu$ of the same sort as before, namely
\begin{equation}\begin{split} &\mu^{(k)} = \nu\rightarrow \cr
 &\mu^{(k)} + 
 \Lambda_{k-1} \rightarrow \dots \rightarrow
\mu^{(k)} +  d_{k-1}' \Lambda_{k-1} = \mu^{(k-1)}  \rightarrow \cr
 &\mu^{(k-1)} + 
 \Lambda_{k-2} \rightarrow \dots \rightarrow
\mu^{(k-1)} +  d_{k-2}' \Lambda_{k-2} = \mu^{(k-2)}  \rightarrow \cr &\vdots
\cr &	\mu^{(2)} + \Lambda_{1} \rightarrow \dots \rightarrow \mu^{(2)} + d_1'
\Lambda_{1} =
\mu.
\end{split}\end{equation} This skew tableau has property (L), so we can
define $\tilde A(\mu) = f(T) \tilde G(\nu)$.

\medskip
\noindent{\bf Case 3. Non-interior diagrams on a critical face. }
\medskip

We continue to use the notation  of the previous case. Next consider the
case that $\mu$ is a non-interior Young diagram, that $\mu + \rho$ is not
contained in the fundamental box $B$, but 
$\mu+ \rho$ already lies in the interior $F^0$ of the face
$F$ of $a + A^+$.  Put $J = \{1,2,\dots, k-1\} \setminus I$.  Note that
$a$ lies on the face $H$ of the Weyl chamber 
$\mathcal C$  generated by $\{\Lambda_j : j
\in J\}$. If there is some $j \in J$ such that $\mu - l \Lambda_j$ is in the
interior of the Weyl chamber, then take the greatest such $j$ and put $\nu =
\mu - l \Lambda_j$. 
(Thus, in going from $\mu + \rho$ to $\nu + \rho$, one moves
toward the origin, parallel to the face $H$ of $\mathcal C$, and
$\nu + \rho$ lies on an opn face of the same type as $F^0$.)
Furthermore, take $T$ to be the skew tableau of shape
$\mu \setminus \nu$
$$
\nu \rightarrow \nu + \Lambda_j \rightarrow \cdots \rightarrow \nu + l
\Lambda_j = \mu,
$$ where at each stage, the j nodes of $\Lambda_j$ are added from top to
bottom. Then $T$ has property (L), and we put $\tilde A(\mu) = f(T)\tilde
G(\nu)$.

Finally, if there is no $j \in J$ such that 
$\mu - l \Lambda_j$ is in the interior of the Weyl chamber, then compute
$\tilde A(\mu)$ by the original LLT algorithm; that is, take $T$ to be the
standard ladder tableau of shape $\mu$ and put
$\tilde A(\mu) = f(T) \emptyset$.

\medskip
\noindent{\bf Case 4. $k > l$ }
\medskip

When $l < k$, only faces of dimension $\le l-1$ contain dominant integral
weights, and not all dominant integral weights are $l$-regular.  No
modification to the algorithm is necessary in the interior region, but the
algorithm for the boundary regions must be modified as follows.  In case 2,
compute
$\nu$ as before, and take the first $\nu' = \mu{k-j} + s \Lambda_{k-j+1}$ on
the canonical tableau from $\nu$ to $\mu$ which is $l$-regular.  Take
$T'$ to be the tail of this tableau from $\nu'$ to $\mu$ and put
$\tilde A(\mu) = f(T') \tilde G(\nu')$.  In case $\mu = \nu'$, proceed as in 
case 3. This completes the description of the algorithm.

\bigskip

Empirically, the modified algorithm produces enormous improvements in
efficiency. To demonstrate this, we include some timing experiments
comparing the original LLT algorithm with the modified algorithm.\footnote{
The algorithms were encoded in Mathematica, and run on a  266 mhz
 Apple Macintosh G3 computer.}  We also compare the algorithm given in
~\cite{Soergel1}; cf. Section 5.\footnote{The Soergel algorithm may have
been  somewhat disadvantaged by  inefficient programming.}  Times are
reported in seconds for computing
$\tilde G(\mu)$ for the given $k$, $l$, and $\mu$. In the table,  $N$
denotes the number of diagrams
$\mu' <
\mu$ for which
$\tilde G(\mu')$ had to be computed in the recursion.

\bigskip
\centerline{
\begin{tabular}{|l|c| c|c| }
\hline
$k,\ l,\ \ \mu$ &LLT,\  secs.,\ N  & modified LLT,\  secs., \ N & Soergel,\ 
secs., \ N 
\\
\hline
\hline
$4 ,5,\  [4 l, 2 l, 0,0]$ & 13.5,\ 16 &1.43,\ 5 &   5.22, \ 18  \\
\hline
$4 ,5, \ [8 l, 4 l, 0,0]$ & 2013,\ 619 &5.95,\ 17  &179,\ 81  \\
\hline
$5, 6, \ [6l, 4 l, 2 l, 0,0]$ &56924,\ 2245 & 8.52,\ 8  & 6771,\ 391 \\
\hline
$5, 6, \ [12l, 8 l, 4 l, 0,0]$& $>$ 24 hrs. & 50,\ 48 & $>$ 24 hrs. \\
\hline
\end{tabular} }
\bigskip

\section{The lower global crystal base and affine Kazhdan-Lusztig
polynomials}

In this section we will demonstrate that the coefficients of the lower global
crystal base for $\usl$ coincide with certain affine Kazhdan-Lusztig
polynomials. We first recall the definition of these Kazhdan-Lusztig
polynomials, following ~\cite{Soergel1}.

Let $R \subseteq \mathfrak h^* \cong \R^{k-1}$  denote the root system of
$\mathfrak{sl}_k$ contained in the dual of the diagonal subalgebra.  Let 
$W
\cong S_k$ denote the Weyl group.  The affine Weyl group is the semi-direct
product
$\W = W \ltimes \Z R$, which acts on $\mathfrak h^*$.  Fix a positive
integer $l$.  We consider the  {\em level $l$} action of the affine Weyl
group on $\mathfrak h^*$, i.e., the action  via the natural isomorphism of 
$\W$ with its subgroup 
$\W^{(l)}  = W \ltimes l \,\Z R$.  $\W^{(l)}$ is generated by reflections in
certain affine hyperplanes in $\mathfrak h^*$.  The connected components of
the complement  of the union of all of these reflection hyperplanes are
called the {\em alcoves} at level $l$.  The set of all alcoves is denoted by 
$\mathcal A$, and the set of all alcoves which are contained in the positive
Weyl chamber $\mathcal C$ is denoted by $\mathcal A^+$.  The alcove
$A^+$ is the unique element of $\mathcal A^+$ which contains the origin
$0$ in its closure.  Let $\S$ denote the set of reflections in the walls of 
$A^+$.  Then $(\W, \S)$ is a Coxeter group with generating set of reflections
$\S$.

The affine Hecke algebra $\H = \H(\W, \S)$ is the associative algebra  with
identity element $\one$ over  the ring of Laurent polynomials
$\Z[v, v\inv]$ with generators $\{T_s : s \in \S\}$ which satisfy the  braid
relations and the quadratic relation $T_s^2 = v^{-2} \one + (v^{-2} -1) T_s$,
for $s \in \S$.  Using instead the generators $H_s = v T_s$, one has instead
the quadratic relations $H_s^2 = H_s + (v - v\inv) \one$.  The Hecke algebra
has a basis $\{H_x : x \in \W\}$ satisfying $H_xH_y = H_{xy}$ in case
$\ell(xy) = \ell(x) + \ell(y)$, and $H_sH_x = H_x + (v-v\inv)H_{sx}$ in case
$\ell(sx) = \ell(x) -1$, for $x \in \W$ and $s \in \S$.  The Hecke algebra
has an involution $d : a \mapsto \bar a$ defined by $\bar v = v\inv$ and
$(H_x)^- = (H_{x\inv})\inv$.  An element fixed by this involution is called
self-dual. Fundamental self-dual elements are the $C_s = H_s + v$.

Let $\S_0 \subseteq \S$ be the set of reflections fixing the origin. The
finite Weyl group $W$ is the Coxeter subgroup of $\W$ with generating set of
reflections $\S_0$. Let $\H_f = \H(W, \S_0) \subseteq \H$ denote its Hecke
algebra. Consider the sign representation of $\H_f \rightarrow \Z[v, v\inv]$
which takes each $H_s$ to $-v\inv$.  Let $\N$ denote the induced right
$\H$-module
$$
\N = \Z[v, v\inv] \otimes_{\H_f} \H
$$  Let $\W^f \subseteq \W$ be the coset  representatives of minimal length
of the right cosets of $W$ in $\W$.  Then $\N$ has a basis
$N_x = \one \otimes H_x$, for $x \in \W^f$, and the operation of  the $C_s$
for $s \in \S$ on this basis has the following form:
\begin{equation} N_x C_s =
\begin{cases} N_{xs} + v N_x &\text{if $xs \in \W^f$ and $xs > x$;}\\ N_{xs}
+ v\inv N_x &\text{if $xs \in \W^f$ and $xs < x$;}\\ 0   &\text{if $xs
\not\in \W^f$}
\end{cases},
\end{equation} where the inequality signs refer to the Bruhat order on
$\W$.  The involution on $\H$ induces an involution on $\N$ defined by $a
\otimes b \mapsto \bar a \otimes \bar b$.

The following result is  discussed in ~\cite{Soergel1}, following
~\cite{Deodhar1},~\cite{Deodhar2}. 

\def\Nbar{\underline N}

\medskip
\vbox{
\begin{theorem}  $\N$ has a unique basis $\{\underline N_x : x \in \W^f\}$
satisfying
\begin{enumerate}
\item $\Nbar_x \in N_x + \sum_{y < x} v\Z[v] N_y$.
\item $\Nbar$ is self-dual.
\end{enumerate} 
\end{theorem} }

\medskip
\begin{definition}  The affine Kazhdan-Lusztig polynomials $n_{y,x}$ are
defined by
$$
\Nbar_x = N_x + \sum_{y < x} n_{y, x} N_y
$$
\end{definition}

\medskip The affine Weyl group $\W$ acts freely and transitively on alcoves,
so there is a bijection $\W \rightarrow \A$ given by $w \mapsto w A^+$, where
$A^+$ is the unique alcove in $\A^+$ containing 0 in its closure. Under this 
bijection, the elements of $\W^f$ correspond to alcoves contained in the 
positive Weyl chamber.  One also has an action of $\W$ on the right on $\A$
given by
$(w A^+) x = wx A^+$.  

Using the bijection between $\W^f$ and $\A^+$, one can  rename the
distinguished elements of the right $\H$ module
$\N$ using alcoves $A \in \A^+$ rather than coset representatives $x \in
\W^f$.   Thus if $x, y \in \W^f$ correspond to $A, B \in \A^+$, then we write
$N_A$ for $N_x$, $\Nbar_A$ for $\Nbar_x$, and $n_{A, B}$ for $n_{x, y}$. The
right action of $\H$ is then given by

\begin{equation} \label{equation right action of Cs} N_A C_s =
\begin{cases} N_{As} + v N_A &\text{if $As \in \A^+$ and $As \succ A$;}\\
N_{As} + v\inv N_x &\text{if $As \in \A^+$ and $As \prec A$;}\\ 0  
&\text{if $As \not\in \A^+$},
\end{cases},
\end{equation} where now the inequalities have a geometric interpretation:
$As \succ A$ if $As$ is on the positive side of the hyperplane separating
the two alcoves.  We remark that the $\Nbar_A$ are computed by a recursive
scheme reminiscent of the computation of the lower global crystal base. One
has $\Nbar_{A^+} = N_{A^+}$.  Given $A \ne A^+$, one can choose $s \in \S$
such that $As \in \A^+$ and $As \prec A$.  As a first approximation to
$\Nbar_A$ one takes
$$
\Nbar_{As} C_s = N_A + \sum_{B \prec A} f_{B, A}(v) N_B.
$$ This element is self-dual, but may have coefficients with non-zero
constant term.  So one corrects these coefficients by subtracting a self-dual
linear combination of $\Nbar_B$ for $B \prec A$.

Finally, we want to rename the elements of $\N$ one last time, using dominant
integral weights of $\mathfrak{sl}_k$.  We fix an integer $l$ and take the
level $l$ action of the affine Weyl group on $\R^k$, generated by
reflections in the hyperplanes
$x_i - x_j = m l$.    Suppose first that a weight $\mu$ lies in an open
alcove.  Write $a(\mu)$ for the alcove of $\mu$ and define
$$ n_{\la, \mu} = 
\begin{cases} n_{a(\la), a(\mu)} &\text{if $\la$ is in the $\W$ orbit of
$\mu$;}\\ 0 & \text{otherwise,}
\end{cases}
$$ and  $\Nbar_\mu = \sum_\la n_{\la, \mu} \la$. Now consider a weight $\mu$
which lies on one or more affine hyperplanes, and define $a^+(\mu)$ to be
the unique open alcove which contains $\mu$ in its closure and which lies on
the positive side of all hyperplanes containing
$\mu$.  In this case we put
$$ n_{\la, \mu} = 
\begin{cases} n_{a^+(\la), a^+(\mu)} &\text{if $\la$ is in the $\W$ orbit of
$\mu$;}\\ 0 & \text{otherwise,}
\end{cases}
$$ and again $\Nbar_\mu = \sum_\la n_{\la, \mu} \la$.  (Note that for these
definitions, it is quite unnecessary to assume $k \le l$ so that the open
alcoves in fact contain integral weights.)

\begin{theorem} Fix integers $k$ and $l$. Let $\la$ and $\mu$ be Young
diagrams of the same size, both with no more than
$k$ rows, and with
$\mu$ $l$-regular.  Then $\dlm(v) = n_{\la + \rho, \mu + \rho}(v)$.
\end{theorem}

The proof of this result, while straightforward, will require a number of
intermediate lemmas and observations.  We begin by recalling a geometric
interpretation of tableaux and of conjugacy of tableaux (with no more than
$k$ rows.)  To each Young diagram $\la$ there corresponds the dominant
integral weight
$\tilde \la$ of $\mathfrak{sl}_k$ given by $\tilde\la_i = \la_i - \la_k$.

To a  standard (skew) tableau 
$$ T = (\la^{(0)} \subseteq \la^{(1)} \subseteq \cdots \subseteq \la^{(s)}),
$$ there corresponds the {\em path}
$$ (\tilde\la^{(0)} \subseteq \tilde\la^{(1)} \subseteq \cdots \subseteq
\tilde\la^{(s)}),
$$ in the positive Weyl chamber; one might picture the succesive $\tilde
\la^{(i)}$ as begin connected by affine line segments, so that the path in
fact becomes a piecewise linear curve.  The tableau may be recovered from
the path and the initial diagram $\la^{(0)}$, so we will not distinguish
between tableaux and paths.  The tableau $T+ \rho$ is that obtained by
adding the half sum of positive roots $\rho$ to each diagram $\la^{(i)}$. 
Note that the length $\la_i -i +k$ of the  $i\th$ row of $\la + \rho$ is the
content of the node $(i, \la_i)$ of $\la$, plus $k$.   Therefore, the weight
of a tableau $T$, computed in Section 3 in terms of contents of nodes of
diagrams along the path $T$, can be computed instead in terms of row lengths
of diagrams along the path  $T+\rho$.  As the hyperplanes of the affine
reflection group $\W$ at level $l$ are the loci of points having two
particular coordinates conjugate modulo $l$, two standard tableaux $S, T$ are
$l$-conjugate if and only if the corresponding paths $S+\rho, T + \rho$ are
related by reflections in such hyperplanes.  In particular, the end-points of
such paths are in the same $\W$ orbit. For our present purposes it will
suffice to consider standard skew tableaux whose residue sequences (cf.
Definition 3.2) have all multiplicities $m_i$  equal to $1$.

In the following, a {\em face} will always mean a face of $\overline A$,
where
$A \in \A^+$ is an open alcove.  An {open face} is the interior of such a
face; cf. Section \ref{section preliminaries}.  Given faces $F$ and $F'$ in a
$\W$-orbit,  write $F' \trianglelefteq F$ if for all $\mu\in F$ and $\mu'
\in 
\W \mu \cap F'$, one has $\mu' \trianglelefteq \mu$.

\medskip
\begin{lemma} \label{lemma dependence on faces only}   Let $T$ be a standard
skew tableau such that the corresponding path $T + \rho$ is contained in
some open face.  Then the weight of
$T$ is
$$
\wt(T) = 1.  
$$ Furthermore, $T$ has no $l$-conjugates (other than itself).
\end{lemma}

\begin{proof}  At the $i\th$ step in the tableau $T$ a node of some residue
$r_i$ is added to a diagram $\la^{(i-1)}$. It follows from the assumption
that $T$ stays in the open face $F$, that $\la_i$ has no other indent nodes
or removable nodes of the same residue.  Therefore the weight of $T$ is $1$
and $T$ has no 
$l$-conjugates other than itself.
\end{proof}

\medskip We will work for a while under the assumption that $k \le l$, so
that open faces of all dimensions contain dominant integral weights and all
Young diagrams of length $\le k$ (except for $n \Lambda_k$) are 
$l$-regular.

\medskip
\vbox{
\begin{corollary}  
\label{corollary dependence on faces only}   Let $\mu, \la, \mu',  \la'$
be Young diagrams,  such that
$\la+\rho$ is in the $\W$ orbit of $\mu+\rho$ and
$\la'+\rho$ is in the $\W$ orbit of $\mu'+\rho$.
If 
$\mu + \rho$ and $\mu' + \rho$ lie on the same open face, 
and 
$\la + \rho$ and $\la' + \rho$ lie on the same open face,
then 
$$\dlm(v) = d_{\la', \mu'}(v).$$ 
\end{corollary} }

\begin{proof} One can assume without loss of generality that 
$\mu' \subseteq \mu$, and that there is a standard skew tableau
$T$ of shape $(\mu + \rho) \setminus (\mu'+ \rho)$ which lies  entirely
contained in the open face containing the two endpoints.   It follows from
Lemma 
\ref{lemma dependence on faces only} that  
$\tilde G(\mu) = f(T)\tilde G(\mu') = \sum_{\la'} d_{\la'
\mu'}(v) f(T)
\la' =
\sum_\la d_{\la' \mu'}(v) \la$, where the first sum is over the $\W$ orbit
of $\mu'$ and the second over the $\W$ orbit of $\mu$.
\end{proof}

\medskip
\begin{definition} Let $F$, $F'$  be  open faces in  the same $\W$ orbit 
such that  $F'
\trianglelefteq F$. Suppose  $\la$ and $\mu$ are Young diagrams of the same
size, $\mu+\rho \in F$, and  $\la +\rho \in F' \cap \W \mu$.   Define
$d_{F', F}(v)$ to be $\dlm(v)$, and
$$
\tilde G(F) = \sum_{F'} d_{F', F} F'.
$$
\end{definition}

\medskip This makes sense according to Corollary \ref{corollary dependence
on faces only}.

\ignore{
\medskip
\begin{definition} Let $F$ be an open face of dimension $d \le k-1$.  If $d
= k-1$, then
$F$ is an open alcove, and  for $\mu \in F$, no two rows of $\mu$ are
conjugate modulo $l$.  If $d < k-1$, then there exist subsets $I_1, I_2,
\dots, I_s$ of
$\{1, 2, \dots, k\}$ of cardinalities $n_1, n_2, \dots, n_s$ such that
\begin{enumerate}
\item $2 \le n_1 \le n_2 \dots \le n_s$,
\item $\sum_i (n_i -1) = k-1-d$, and
\item For $\mu \in F$, $\mu_i \equiv \mu_j \quad (\mod l)$ if and only if
$i, j \in I_r$ for some $r$.
\end{enumerate} We will say that $F$ is an {\em open face of type} $(n_1,
n_2, \dots, n_s)$.
\end{definition} }

\medskip
\begin{lemma} \label{lemma leaving faces} Let $F$ be an open face of
dimension $d < k-1$.  Let 
$I  \subseteq \{1, 2, \dots, k\}$  be a subset of cardinality
$\ge 2$ which is maximal  with respect to the property that  for all $\mu
\in F$
$\mu_i \equiv \mu_j \quad (\mod \ l)$ if 
$i, j \in I$ .  Then there is a $\mu \in F$ such that for all $i \in I$,
$\mu + \eps_i$ lies in an open face of dimension $d+1$.
\end{lemma}

\begin{proof}  Left to the reader.
\end{proof}

\medskip Let $F$ be an open face of dimension $d < k-1$, and let 
$I  \subseteq \{1, 2, \dots, k\}$  be a subset with the property described
in the lemma.  For each
$i
\in I$ let
$F_i$ be the open face of dimension $d+1$ which contains $\mu + t\eps_i$ for
$\mu \in F$ and for $t > 0 $ small.   Set $F_+(I) = F_{i_0}$, where $i_0$ is
the least element of $I$.
$F_+(I)$ has the property:   For each reflection hyperplane $H$ of $\W$
which contains $F$ but which does not contain $F_+(I)$, the face  $F_+(I)$
lies on the positive side of $H$. 

Let  $F'$ be an open face such that $F'$ lies in the $\W$ orbit of $F$ and
$F' \trianglelefteq F$.
 Let  $I'$ be the set of indices which plays the same role for $F'$ as does
$I$ for $F$.  Namely,  if $w \in \W$ is an element such that $w F = F'$, let 
$w(x_1, \dots, x_k) = (x_{\sigma\inv(1)}, \dots, x_{\sigma\inv(k)}) + w(0)$
for all $(x_1, \dots, x_k) \in \R^k$, and put $I' = \sigma(I)$.  Note that
$\sigma$ is not uniquely determined, but 
$I'$ is uniquely determined by $F$ and $F'$, and $I$.  Define
$F'_+(I')$ in the same way as $F_+(I)$, namely 
$F'_+(I') = F'_{j_0}$, where $j_0$ is the least element of $I'$.

\medskip
\begin{lemma}  $$d_{F ', F}(v) = d_{F'_+(I'), F_+(I)}(v).$$
\end{lemma}

\begin{proof} Adopt the notation of the two paragraphs preceding the
statement of the lemma. Let $\mu$ be a Young diagram such that $\mu + \rho
\in F$ and such that
$\mu + \rho + \eps_{i_0} \in F_+(I)$, as is possible by Lemma 
\ref{lemma leaving faces}.  Let $r$ be the residue of the indent node of 
$\mu$ in the row $i_0$.  By the assumption on $\mu$, the Young diagram $\mu$
has exactly $|I|$ indent $r$-nodes and no removable $r$-nodes, and the same
holds for any $\la$ such that $\la + \rho$ is in the $\W$ orbit of $\mu +
\rho$. For any such $\la$, with $\la + \rho$ lying on an open face $F'$, one
has
$$ f_r \la = \sum_{i \in I'} v^{N(\la, \la + \eps_i)} (\la + \eps_i) ,
$$ and furthermore $N(\la, \la + \eps_i) \ge 0  $, with equality if and only
if
$\la + \eps_i \in F'_+(I')$.  It follows that
$$ f_r \tilde G(\mu) = \tilde G(\mu + \eps_{i_0}),
$$ and furthermore
$$ d_{F'_+(I'), F_+(I)}(v) = d_{\la + \eps_{j_0}, \mu +\eps_{i_0}}(v)  =
\dlm(v),   = d_{F', F}(v).
$$
\end{proof}

\medskip For a face $F$, let $a^+(F)$ denote the unique alcove $A$ such that 
$F$ is contained in the closure of $A$ and $A$ lies on the positive side of
all hyperplanes containing $F$.

\medskip
\begin{corollary} For any face $F$ and any face $F'$ in the $\W$ orbit of
$F$ such that $F' \trianglelefteq F$, one has
$$ d_{F ', F}(v) = d_{a^+(F'), a^+(F)}(v)
$$
\end{corollary}

\begin{proof}  This follows by induction from the previous lemma.
\end{proof}

\medskip To complete the proof of Theorem 5.3 in the case $k \le l$, it
remains to show that $d_{B, A} = n_{B, A}$ for alcoves $B, A$.   The $n_{B,
A}$ are defined by a recursion involving crossing of walls, so we will have
to see that the $d_{B, A}$ satisfy the same recursion.

By a wall of an alcove, we mean a face of the alcove of dimension $k-2$,
that is the non-empty intersection of the closure of the alcove with a
reflection hyperplane of $\W$.

\medskip
\vbox{
\begin{lemma} \label{lemma wall crossing} Let $A$ and $B$ be adjacent
alcoves, separated by an open wall
$F$, and let $T$ be a skew tableau such that $T+\rho$ which starts at a
diagram $\nu+\rho \in A$ and ends in a diagram $\mu + \rho \in B$, passing
through $F$. Denote the reflection of $\mu + \rho$ in $F$ by  $\mu' + \rho$.
\begin{enumerate}
\item If $A \prec B$ then $f(T) \nu = \mu + v \mu'$.
\item If $A \succ B$ then $f(T) \nu = \mu + v \inv\mu'$.
\end{enumerate}
\end{lemma} }

\begin{proof} $T$ has one $l$-conjugate tableau $T'$ which ends in $\mu'$.
It has to be shown that $T$ has weight $\wt(T) = 1$, while
$T'$ has weight $\wt(T') = v$ in case (a), and
$\wt(T') = v\inv$ in case (b).
To verify this, one has to check what happens when a 
path hits a wall from above or below, or when a path leaves
a wall towards the positive or negative side.  The various
cases to be checked are listed in the following lemma.
\end{proof}

\medskip
\begin{lemma}  Let $T$ be a skew tableau.
\begin{enumerate}
\item  If $T+\rho$ begins in an alcove $A$ and ends in an open wall $F$ of
the alcove, and $A$ lies on the negative side of the hyperplane containing
$F$, then $\wt(T) = 1$.
\item If $T+\rho$ begins in an alcove $A$ and ends in an open wall $F$ of
the alcove, and $A$ lies on the positive side of the hyperplane containing
$F$, then $\wt(T) = v\inv$.
\item If $T+\rho$ begins in an open wall $F$ of an alcove $A$ and ends in
the alcove $A$, and  $A$ lies on the negative side of the hyperplane
containing $F$, then $\wt(T) = v$.
\item If $T+\rho$ begins in an open wall $F$ of an alcove $A$ and ends in
the alcove $A$, and  $A$ lies on the positive side of the hyperplane
containing $F$, then $\wt(T) = 1$.
\ignore{
\item  If $T+\rho$ begins in an alcove and ends in an adjacent alcove, 
passing through the open wall dividing the two alcoves, then
$\wt(T) = 1$.
\item If $T+\rho$ begins and ends in an alcove $A$, but touches an open wall
$F$ of the alcove, and $A$ lies on the negative side of the hyperplane
containing $F$, then $\wt(T) = v$.
\item If $T+\rho$ begins and ends in an alcove $A$, but touches an open wall
$F$ of the alcove, and $A$ lies on the positive side of the hyperplane
containing $F$, then $\wt(T) = v\inv$.
}
\end{enumerate}
\end{lemma}

\begin{proof} \ignore{Parts (e)-(g) follow from the previous parts, together
with  Lemma \ref{lemma dependence on faces only}.}  The proofs of parts
(a)-(d) are all similar, so   we prove part (b) and leave the rest to the
reader.   Let
$$ T = (\la^{(0)} \subseteq \la^{(1)} \subseteq \cdots \subseteq \la^{(s)}),
$$ where $\la^{(i)} + \rho \in A$ for $i < s$ and $\la^{(s)} + \rho \in F$.
It follows from Lemma \ref{lemma dependence on faces only} that
$\wt(T) = v^{N(\la^{(s-1)},\la^{(s)})}$.  Say the wall $F$ is given by
$x_a - x_b = m l$, where $a < b$.  Then $\la^{(s)}$ is obtained from
$\la^{(s-1)}$ by filling an indent node of some residue $r$ in row $b$
(since $A$ lies above $F$), and furthermore
$\la^{(s-1)}$ has no other indent $r$-nodes, but has a removable $r$-node in
row $a$.  It follows that $N(\la^{(s-1)},\la^{(s)}) = -1$ and 
$\wt(T) = v\inv$.
\end{proof}

\medskip
\begin{lemma} For all alcoves $A, B$, 
$$ d_{B, A} = n_{B, A}
$$
\end{lemma}

\begin{proof}  Define a map $\Phi : \mathcal F_k \rightarrow \N$ by  putting
$\Phi(\mu) = N_A$ if $\mu + \rho$ is in the open alcove $A$, and
$\Phi(\mu) = 0$  if $\mu$ is not contained in an open alcove.  The assertion
of the lemma is equivalent to
$\Phi(\tilde G(\mu)) = \Nbar_A$, if $\mu + \rho \in A$.  If $\mu + \rho \in
A^+$, then $\Phi(\tilde G(\mu)) = N_{A^+} = \Nbar_{A^+}$.  The proof
proceeds by induction on $\mathcal A^+$ with respect to $\prec$.  For an
alcove
$A \ne A^+$, choose $s \in \S$ such that $As \prec A$, and let $F$ be the
open wall separating $As$ and $A$.  One can choose a skew tableau $T$ of
shape
$\mu \setminus \nu$ such that $\nu + \rho \in As$, and $\mu + \rho \in A$, 
and $T+\rho$ crosses from $As$ to $A$ through $F$.  Put
$\tilde A(\mu) = f(T) \tilde G(\nu)$.  One has
$$
\Phi(\tilde A(\mu)) = \Phi(\tilde G(\nu)) C_s = \Nbar_{As} C_s,
$$ where the first equality comes from Lemma \ref{lemma wall crossing}, and
the definition of the right action of $C_s$ (Equation 
\ref{equation right action of Cs}), and the second equality  from the
induction hypothesis.  In particular
 $\tilde A(\mu) = \sum_\la
\alm \la$, where $\alm \in \nat[v]$.  Therefore the rectification of $\tilde
A(\mu)$ to $\tilde G(\mu)$ takes the simple form
$$
\tilde G(\mu) = \tilde A(\mu) - \sum_{\la < \mu} \alm(0) \tilde G(\la).
$$ Thus
$$
\Phi(\tilde G(\mu) ) = \Nbar_{As} C_s - \sum_{\la < \mu} \alm(0)
\Nbar_{a(\la+\rho)} = \Nbar_A,
$$ using the induction hypothesis.
\end{proof}

\def\fac{\mathbf F} This lemma completes the proof of Theorem 5.3, in case
$k \le l$.  To handle the case $l < k$, we show that the polynomials
$\dlm(v)$ are independent of $l$ in a certain sense.

Let $\fac^l$ denote the set of all open faces for the action of the affine
Weyl group $\W^{(l)}$ at level $l$.  
 If $l_1 < k <  l_2$, then the map $x \mapsto  l_2 l_1\inv x$ from $\mathcal
C$ to $\mathcal C$ induces a bijection $\psi :
\fac^{l_1} \rightarrow \fac^{l_2}$.  For $\mu \in F \in \fac^{l_1}$ and $\nu
\in \psi(F)$, define a map
$\psi_\mu^\nu$ from the set of $\mu'$ such that $\mu' + \rho$ is in the
$\W^{(l_1)}$ orbit of $\mu + \rho$ to the set of $\nu'$ such that $\nu' +
\rho$ is in the
$\W^{(l_2)}$ orbit of $\nu + \rho$ by $\psi_\mu^\nu (\mu') = \nu'$ if $\mu'
+ \rho \in H \in \fac^{l_1}$ and $\nu' + \rho\in \psi(H)$.  Extend
$\psi_\mu^\nu$ linearly to the $\Z[v, v\inv]$- modules spanned by such
diagrams.

\ignore{ Extend
$\psi$ linearly to the $\Z[v, v\inv]$ modules with bases $\fac^{l_1}$ and 
$\fac^{l_2}$.   Say an open  face $F \in \fac^{l_1}$ is $l_1$-regular if it
contains $\mu + \rho$ for some $l_1$-regular $\mu$. } 

\ignore{
\medskip
\begin{proposition} \label{proposition independence of l} 
$\psi(\tilde G(F)) = \tilde G(\psi(F))$ for all
$l_1$-regular open faces $F$.  
\end{proposition} }

\def\barmu{{\bar \mu}}
\def\barla{{\bar \la}}

\begin{lemma}\label{lemma independence of l}
 Let $\mu$ be an $l_1$-regular diagram,  and $F$ the open face containing
$\mu + \rho$.
 Let $T$ be the standard ladder tableau
 with shape $\mu$.  Put $\tilde A(\mu) = \tilde A(T)$.   There exists a
$\bar \mu$ with $\bar \mu + \rho \in \psi(F)$ such that
$\psi_\mu^\barmu(\tilde A(\mu))$ is self-dual.
\end{lemma}

\begin{proof}  The proof goes by induction on the number of nodes of $\mu$.
If $\mu$ is the empty diagram, let $F$ the the open face containing $\mu +
\rho$ and let $\barmu$ be any diagram such that  $\barmu + \rho \in
\psi(F)$.  Then 
$\psi_\mu^\barmu(\tilde A(\mu)) = \barmu$.  But, since $\barmu$ is the
lowest diagram in its $\W^{(l_2)}$ orbit, $ \barmu = \tilde G(\barmu) $. 
Thus 
$\psi_\mu^\barmu(\tilde A(\mu))$ is self-dual.

Now fix an $l_1$-regular $\mu \in F \in \fac^{l_1}$,  and assume the
assertion holds for all $l_1$-regular diagrams with fewer nodes.  Let $T'$
be the tableau obtained by removing the last node of $T$, $\mu'$ the shape
of $T'$, and $F'$ the open face containing
$\mu' + \rho$.   Note that $\mu'$ is $l_1$-regular. By the induction
hypothesis, there is a diagram $\barmu'$ with $\barmu' + \rho \in \psi(F')$
such that $\psi_{\mu'}^{\barmu' } (\tilde A(T'))$ is self-dual.

Let $r$ be the residue of the node $\mu\setminus\mu'$ and let 
$c = v^{N(\mu', \mu)}$.  Then $$\tilde A(\mu) = c\inv f_r \tilde A(T').$$

We assert the existence of  a skew tableau $\bar S$ such that $\bar S +
\rho$ starts at $\barmu + \rho \in \psi(F')$ and ends in some diagram 
$\barmu + \rho \in \psi(F)$

If $F' = F$, then there is nothing to show.  Otherwise, one might have
$F$ an open face of some dimension $d$ and $F'$ an open face of dimension
$d-1$ contained  in the boundary of the closure of $F$, or $F'$ a face of
dimension $d+1$  such that $F$ is contained in the boundary of its closure. 
In either case,
$\psi(F')$ and $\psi(F)$ stand in the same relation, and a path of the
desired type exists.  Or one might have $F$ and $F'$ both $d$ dimensional
boundary faces  of a
$d+1$ dimensional face $H$.  In this case, a path from 
$\barmu' + 
\rho \in \psi(F') $ to $\psi(F)$ might necessarily contain integral points of
$\psi(H)$,  while the one step path from
$\mu' + \rho$ to $\mu + \rho$ has no points in $H$; this difference,
however, has no effect.

In all cases, one has for all $\la'$ such that $\la' + \rho$ is in the
$\W^{(l_1)}$ orbit of $\mu' + \rho$,
$$
\psi_\mu^\barmu(f_r \la')   = f(\bar S) \psi_{\mu'}^{\barmu'} (\la'),
$$ because the effect of $f(\bar S)$ is determined entirely by the faces
which $\bar S + \rho$ enters and leaves.

Consequently, one has
$$
\psi_\mu^\barmu(\tilde A(\mu)) = c\inv f(\bar S) \psi_{\mu'}^{\barmu' }
(\tilde A(T')).
$$ The desired conclusion follows from this.
\end{proof}

\medskip
\begin{proposition} \label{proposition independence of l} Let $\mu$ be an
$l_1$ regular diagram and let 
$F$ be the open face containing $\mu + \rho$.  For every diagram
$\nu$ such that $\nu  + \rho \in \psi(F)$, 
$
\psi_\mu^\nu(\tilde G(\mu)) = \tilde G(\nu)
$.
\end{proposition}

\begin{proof} Given $\mu$, let $\tilde A(\mu)$ and  $\bar \mu$  be as in the
lemma, and put
$\tilde A(\bar \mu) = \psi_\mu^\barmu(\tilde A(\mu))$. It follows from the
lemma that
$\tilde A(\bar \mu)$ satisfies properties ($\tilde {\rm A}$1) and ($\tilde
{\rm A}$2).

If $\mu + \rho$ is the lexicographically smallest diagram in its
$\W^{(l_1)}$ orbit such that $\mu$ is $l_1$-regular, then
$\tilde A(\mu) = \tilde G(\mu)$.
 It follows that $\tilde A(\bar \mu) \equiv \barmu \quad (\mod \  L_k)$, so
$\tilde A(\bar \mu) = \tilde G(\barmu)$, by the uniqueness of $ \tilde
G(\barmu) $.  Therefore, for such $\mu$, we  have $\psi_\mu^\barmu(\tilde
G(\mu)) = \tilde G(\barmu)$.  But then,  Corollary \ref{corollary dependence
on faces only} implies 
$
\psi_\mu^\nu(\tilde G(\mu)) = \tilde G(\nu)
$ for all $\nu$ such that $\nu + \rho \in \psi(F)$.

Now fix $\mu$ and assume that the result holds for all $\mu' < \mu$ such
that $\mu' + \rho$ is in the $\W^{(l_1)}$ orbit of $\mu + \rho$.  Write
\begin{equation} \label{equation rectification}
\tilde G(\mu) = \tilde A(\mu) - \sum_i \gamma_i(v) \tilde G(\mu^{(i)}), 
\end{equation} where $\mu^{(i)} < \mu$ and $\mu^{(i)} + \rho$ is in the 
$\W^{(l_1)}$ orbit of $\mu + \rho$.  Let $F_i$ be the open face containing
$\mu^{(i)} + \rho$, and let $\nu^{(i)}$ be the diagram with
$\nu^{(i)} + \rho \in \psi(F_i)$ and $\nu^{(i)} + \rho$ in the
$\W^{(l_2)}$ orbit of $\barmu + \rho$.  Then one has
$$
\psi_\mu^\barmu(\tilde G(\mu^{(i)})) = 
\psi_{\mu^{(i)}}^{\nu^{(i)}}(\tilde G(\mu^{(i)})) =
\tilde G(\nu^{(i)}),
$$ with the last equality coming from the induction hypothesis. Applying
this to Equation \ref{equation rectification} gives
$$
\psi_\mu^\barmu(\tilde G(\mu)) = \tilde A(\barmu) - \sum_i \gamma_i(v) \tilde
G(\nu^{(i)}) = \tilde G(\barmu). 
$$ Now applying Corollary \ref{corollary dependence on faces only}   again
gives 
$
\psi_\mu^\nu(\tilde G(\mu)) = \tilde G(\nu)
$
 for all $\nu$ such that
$\nu + \rho \in \psi(F)$.
\end{proof}

This proposition completes the proof of  Theorem 5.3.

\end{document}